\documentclass[11pt]{article}   
\usepackage{amssymb,amscd,latexsym}   
\usepackage{amsmath}
\usepackage{amsthm}
\usepackage{mathdots}
\usepackage[pagebackref]{hyperref}
\usepackage{color}
\usepackage[all]{xy}
\newcommand{\lar}{\longrightarrow}
\newcommand{\llar}{-\kern-5pt-\kern-5pt\longrightarrow}
\newcommand{\surjects}{\twoheadrightarrow}

\newtheorem{Theorem}{Theorem}[section]
\newtheorem{Lemma}[Theorem]{Lemma}
\newtheorem{Corollary}[Theorem]{Corollary}
\newtheorem{Proposition}[Theorem]{Proposition}
\newtheorem{Remark}[Theorem]{Remark}
\newtheorem{Example}[Theorem]{Example}

\newtheorem{Definition}[Theorem]{Definition}
\newtheorem{Question}[Theorem]{Question}

\def\sqr#1#2{{\vcenter{\hrule height.#2pt
        \hbox{\vrule width.#2pt height#1pt \kern#1pt
            \vrule width.#2pt}
        \hrule height.#2pt}}}
\def\phi{\varphi}
\def\demo{\noindent{\bf Proof. }}
\def\square{\mathchoice\sqr64\sqr64\sqr{4}3\sqr{3}3}
\def\qed{\hspace*{\fill} $\square$}

\DeclareMathOperator{\Spec}{Spec}

\DeclareMathOperator{\depth}{depth}

\DeclareMathOperator{\Ht}{ht}
\DeclareMathOperator{\reg}{reg}

\DeclareMathOperator{\grade}{grade}

\def\yy{{\bf y}}

\def\ZZ{{\bf Z}}

\def\fm{{\mathfrak m}}

\def\Ht{{\rm ht}\,}
\def\depth{{\rm depth}\,}
\def\ass{{\rm Ass}\,}

\def\edim{{\rm edim}\,}

\def\grade{{\rm grade}\,}

\def\restr{{\kern-1pt\restriction\kern-1pt}}

\def\pp{{\mathbb P}}

\begin{document}
\begin{center}
{\Large{\bf\sc Equigenerated ideals of analytic deviation one}}
	\footnotetext{AMS Mathematics
		Subject Classification (2010   Revision). Primary 13A02, 13A30, 13D02, 13H10, 13H15; Secondary  14E05, 14M07, 14M10,  14M12.} 
	\footnotetext{	{\em Key Words and Phrases}: plane reduced points, special fiber, perfect ideal of height two, Rees algebra, associated graded ring,  Cohen--Macaulay.}

	\vspace{0.3in}
	
	{\large\sc Zaqueu Ramos}\footnote{Under a post-doc fellowship from INCTMAT/Brazil (88887.373066/2019-00)} \quad\quad
	{\large\sc Aron  Simis}\footnote{Partially
		supported by a CNPq grant (302298/2014-2).}

\end{center}


\begin{abstract}
	\noindent	
The overall goal is to approach the Cohen--Macaulay property of the special fiber $\mathcal{F}(I)$ of an equigenerated homogeneous ideal $I$ in a standard graded ring over an infinite field.
When the ground ring is assumed to be local, the subject has been extensively looked at. Here, with a focus on the graded situation, one introduces two technical conditions, called respectively, {\em analytical tightness} and {\em analytical adjustment}, in order to approach the Cohen--Macaulayness  of $\mathcal{F}(I)$. A degree of success is obtained in the case where $I$ in addition has analytic deviation one, a situation looked at by several authors, being essentially the only interesting one in dimension three.
Naturally, the paper has some applications in this case.	
	
\end{abstract}

\section*{Introduction}

Let $R$ denote a standard graded algebra over a field.
The main focus in this work is on the behavior of the special fiber $\mathcal{F}(I)$ of an equigenerated homogeneous ideal $I\subset R$ of analytic deviation one, in its tight relation with  the associated graded ring ${\rm gr}_I(R)$ and, quite often, the Rees algebra $\mathcal{R}_R(I)$. 
A main target is the Cohen--Macaulayness of $\mathcal{F}(I)$, a venerable subject studied by many authors, specially in the case when the ground ring is local.
Here one approaches the problem presuming that the equigeneration fact ought to yield additional output, a deed of a more recent vintage among experts.

Throughout this introduction, let $I\subset R$ be a $d$-equigenerated ideal, i.e., $I=([I]_d)$.
Two technical conditions are central in this work, as related to the Cohen--Macaulayness of $\mathcal{F}(I)$. The first is the notion of an {\em analytically tight} sequence in power $n$, namely, letting $\ell(I)$ denote the analytic spread of $I$, this is a sequence of $\ell(I)$ forms in $I_d$ having a regularity behavior with respect to powers of $I$.  Traces of the nature of such sequences are spread out in a huge number of important papers -- see \cite{AHT1995}, \cite{CorZar}, \cite{EU2012}, \cite{HubHun}, \cite{HucHun}, \cite{RossiTrungs}, \cite{Trung1987}, \cite{Trung1998}, among uncountable many others.
We make special note of \cite{Trung1998} for the close relationship of its notion of sequences of regular type and the present notion of analytically tight sequences. However, though the technology becomes similar in various passages, the overall objectives are quite different, since here the focus is on the special fiber. By and large the present notion is light enough to handle and does not involve introducing a priori superficial sequences, as we believe is unnecessary for understanding the special fiber of an equigenerated ideal.
The idea of using this concept has a considerable degree of success in the case of an ideal of analytic deviation one in that its existence is guaranteed asymptotically along the powers provided the ideal is unmixed. 
Now, a good deal of previous work assumes as standing hypothesis that the given ideal is generically a complete intersection.
Under such a hypothesis, equigenerated ideals of analytic spread one admit analytically tight sequences in any power, a phenomenon that may explain the success of this popular generic constraint.

Perhaps a true novelty is the use of a second notion, that of an {\em adjusted} set $L$  of forms in $I_d$ and, more particularly, that of an {\em analytically adjusted} such set. The concept is stated in terms of the number of generators of $LI$  and at first sight looks simplistic. Nevertheless it has the merit of fooling quadratic relations. The latter, as known by the experts, is one of the nearly insurmountable obstacles in the theory, usually dribbled by assuming a $G_s$ type of condition.
As the first notion above, it has a strong bearing to properties of minimal reductions of $I$.

Instead of boring the reader with a tedious long introduction, we now briefly describe the contents of each section.

The first section is a quick recap of the main algebras to be used in the work, with a reminder of the special behavior of the special algebra in the equigenerated case.
Here, the first few batches  concern general aspects of any equigenerated homogeneous ideal, so it can help tracking steps in directions other than the one undertaken here.
It includes a few additional aspects of the so-called condition $G_s$ of Artin--Nagata in order to stress some slightly bypassed behavior.
The first of these is a more precise reformulation, in terms of this condition, of the Valla dimension of the symmetric algebra -- a result essentially obtained in the work of Huneke--Rossi (\cite{HR}).
We give some easy consequences of this reformulation.
A seemingly overshadowed result concerns a lower bound for the initial defining degree of the fiber $\mathcal{F}(I)$ of a height two perfect ideal in terms of the $G_s$ condition, a result essentially due to Tchernev (\cite{Tchernev}).

The second section is the core of the paper. It focuses on properties of $\mathcal{F}(I)$ as related to  ${\rm gr}_I(R)$, in particular when the latter is Cohen--Macaulay.
Thus, one proves that, for an equigenerated ideal $I$ of analytic deviation one, if ${\rm gr}_I(R)$ is Cohen--Macaulay then the same property for $\mathcal{F}(I)$ is dictated by the existence of analytically tight sequences in any power. This result may be previously known in the case where $I$ is moreover generically a complete intersection, but we could not find it otherwise in this exact form.
As to the notion of analytical adjustment, it is naturally related to properties of the special fiber $\mathcal{F}(I)$.
Thus, for example, a subset of $I_d$ which is regular sequence over $\mathcal{F}(I)$ is adjusted. Another important example is that of any linearly independent subset of $I_d$ when $I$ is a syzygetic ideal or a perfect ideal of height two  satisfying condition $G_3$.
A basic natural appearance of an analytically adjusted set is as minimal generators of a minimal reduction under certain conditions. 
By and large, if $\mathcal{F}(I)$ is Cohen--Macaulay, any minimal set of generators of a minimal reduction $J\subset I$ is analytically adjusted.
This is a consequence of the freeness of the corresponding Noether normalization.
It is natural to search for conditions under which a converse statement holds. At the moment we could only obtain such a result in dimension three, treated in the last section.

When $\dim R=3$, some of the results become particularly enhanced.
This is the content of the last part of the paper.
 Here both analytical tightness and analytical adjustment fall in one and the same basket to characterize Cohen--Macaulayness of $\mathcal{F}(I)$. The section is  also focused on the reduction number of the ideal $I$ and the multiplicity of its special fiber, with a side about the rational map defined by the linear system of $I$ in its generating degree.


\section{Sprouts of the $G_s$ condition}

Let $R$ be a Noetherian ring and let $I\subset R$ denote an ideal of $R$.

We consider in this work the following graded algebras associated to the pair $(R,I):$

\begin{enumerate}
	\item[$\bullet$] The Rees algebra $\mathcal{R}_R(I)=\bigoplus_{i\geq 0}I^it^i\subset R[t]$ of $I.$ 
		\item[$\bullet$] The associated graded ring ${\rm gr}_I(R)= \mathcal{R}_R(I)/I\mathcal{R}_R(I)$  of $I.$
\end{enumerate}

These algebras are related by a natural $R$-algebra surjection of graded  $R$-algebras $$\mathcal{R}_R(I)\surjects {\rm gr}_I(R).$$
This allows comparison in many directions.
Thus, for  example, in many situations their dimensions differ by one.

In this work we assume throughout that $(R,\fm)$ is a standard graded algebra over a field $k$ with maximal irrelevant ideal $\fm=([R]_1)$.
In this case, if $I\subset \fm$ is a homogeneous ideal then $\dim {\rm gr}_I(R)=\dim R$, while $\dim \mathcal{R}_R(I)=\dim R+1$ provided $\grade I\geq 1$.

In addition, one defines the {\em special fiber cone} (or the {\em  fiber cone}) of $I$ to be 
$$
\mathcal{F}_R(I): = \mathcal{R}_R(I)/\fm\mathcal{R}_R(I),
$$
and the {\em analytic spread} of $I$, denoted by $\ell(I)$, to be the (Krull) dimension of $\mathcal{F}_R(I)$.

Let $I$ be a homogeneous ideal of $R$ minimally generated by forms $f_1,\ldots,f_r$ of the same degree $d$ -- one then says that $I$ is {\em $d$-equigenerated} for short. 

Consider the bigraded $k $-algebra $R \otimes_k  k[y_1,\ldots,y_r]=R[y_1,\ldots,y_r],$ 
where $k[y_1,\ldots,y_r]$ is a polynomial ring over $k$ and the bigrading is given by  ${\rm bideg}([R]_1)=(1,0)$ and ${\rm bideg}(y_i)=(0,1)$.
Setting ${\rm bideg}(t)=(-d, 1)$, then $\mathcal{R}_R(I)=R[It]$ inherits a bigraded structure over $k $.
One has a bihomogeneous $R$-homomorphism
$S \longrightarrow  \mathcal{R}_R(I) \subset R[t] ,
\quad y_i  \mapsto  f_it.$
Thus, the bigraded structure of $\mathcal{R}_R(I)$ is given by 
$$
\mathcal{R}_R(I) = \bigoplus_{c, n \in \ZZ} {\left[\mathcal{R}_R(I)\right]}_{c, n} \quad \text{ and } \quad {\left[\mathcal{R}_R(I)\right]}_{c, n} = {\left[I^n\right]}_{c + nd}t^n.
$$
One is often interested in the $R$-grading of the Rees algebra. Thus, one sets, namely,
${\left[\mathcal{R}_R(I)\right]}_{c}=\bigoplus_{n=0}^\infty {\left[\mathcal{R}_R(I)\right]}_{c,n}$, and of particular interest is  
$$
{\left[\mathcal{R}_R(I)\right]}_{0}=\bigoplus_{n=0}^\infty {[I^n]}_{nd}u^n = k [[I]_d u] \simeq k \left[[I]_d\right]=\bigoplus_{n=0}^\infty {\left[I^n\right]}_{nd}\subset R.
$$
Clearly, $\mathcal{R}_R(I)=[\mathcal{R}_R(I)]_0\oplus \left(\bigoplus_{c\geq 1} [\mathcal{R}_R(I)]_c\right)=[\mathcal{R}_R(I)]_0\oplus \fm \mathcal{R}_R(I)$. Therefore, one gets
\begin{equation}\label{fiber_is_subalgebra}
k \left[[I]_d\right]\simeq {\left[\mathcal{R}_R(I)\right]}_0\simeq \mathcal{R}_R(I)/\fm \mathcal{R}_R(I)=\mathcal{F}_R(I)
\end{equation}
as graded $k $-algebras.

In the sequel, for simplicity, we often write $I_d=[I]_d$ if no confusion arises.


Consider in addition a free $R$-module presentation of $I$:
\begin{equation}\label{syzI}
R^{m}\stackrel{\phi}\lar R^r\to I\to 0.
\end{equation} 

The following notion was introduced in \cite[Section 2]{AN} (see also \cite[Definiton 1.3]{HU}):

\begin{Definition}\rm
	Let $R$ be a Noetherian ring and let $I\subset R$ be an ideal.
Given an integer $s\geq 1$, one says that $I$ satisfies condition $G_s$ if  $\mu(I_p)\leq \Ht p$ for every prime $p\in V(I)$ with $\Ht p\leq s-1.$
\end{Definition}

\begin{Remark}\rm
As is known, this sort of condition and its analogues have an impact on certain dimension theoretic aspects. It is more effective under additional hypotheses, such as when $R$ is Cohen--Macaulay and/or when $I$ has a rank, in which case its effectiveness is expressed in terms of the ideals of minors of a matrix $\phi$ of a presentation of $I$ as in (\ref{syzI}).
\end{Remark}

Throughout this part, $R$ will denote a Noetherian ring of finite Krull dimension $\mathfrak{d}$.
In this context, the most notable cases of the above condition have $s=\mathfrak{d}$ and $s=\mathfrak{d}+1$, the second also denoted $G_{\infty}$.
Here we intend to consider lower cases of $s$, to see how they impact in a few situations.

The symmetric algebra above is very sensitive to it (see, e.g., \cite[Section 7.2.3 and ff]{SiBook}).
The next result is an encore of \cite[Corollary 2.13]{HR} in the present terminology.
The conclusion of the statement was originally conceived by G. Valla, which spurred the subsequent work by C. Huneke and M. E. Rossi.

For the present purpose, recall that a Noetherian ring of finite Krull dimension is said to be {\em equicodimensional} if all its maximal ideals have the same height ($=\dim R$).

\begin{Proposition}\label{Valla's bound} {\rm (The Valla dimension)}
Let $R$ be a Noetherian ring of finite Krull dimension $\mathfrak{d}$ and let $I\subset R$ be an ideal generated by $r$ elements.
Suppose that the following conditions hold:
\begin{enumerate}
	\item [{\rm (i)}] $R$ is equicodimensional and $\grade I\geq 1$.
	\item [{\rm (ii)}] The number of generators of $I$ locally in height $\mathfrak{d}-1$ is at most $r-1$.
	\item [{\rm (iii)}] $I$ satisfies condition $G_{\mathfrak{d}-1}$.
\end{enumerate}
Then $\dim \mathcal{S}_R(I)=\max \{\mathfrak{d}+1,r\}.$
\end{Proposition}
\demo The proof, based on the Huneke--Rossi formula, follows the same pattern with the more or less obvious adaptation.
Thus, first since $\grade I\geq 1$ then for a minimal prime $P\in \Spec R$ of maximal dimension one has $I\not\subset P$, and hence $\dim R/P+\mu(I_P)=\dim R+1=\mathfrak{d}+1$, so  $\dim \mathcal{S}_R(I)\geq \mathfrak{d}+1$.
On the other hand, $\dim \mathcal{S}_R(I)\geq r$ (\cite[Corollary 2.8]{HR}).

Therefore, it suffices to show the upper bound $\dim \mathcal{S}_R(I)\leq \max \{\mathfrak{d}+1,r\}.$
For $I\not\subset P$ there is nothing to prove, so assume that $I\subset P$.
If $P$ is a maximal ideal, equicodimensionality takes care since $\dim R/P=0$ and $\mu(I_P)\leq r$.
If $P$ has height $\mathfrak{d}-1$ assumption (ii) takes care again since $\dim R/P\leq 1$.
Finally, for $P$ of height at most $\mathfrak{d}-2$ the $G_{\mathfrak{d}-1}$ assumption takes over.
\qed

\begin{Remark}\rm
It is curious that a small piece of the proof in \cite[Corollary 2.13]{HR} depends on knowing that the Noetherian base ring $R$ is equicodimensional. This hypothesis is understated in \cite[Corollary 2.13]{HR} possibly because the authors assume throughout that $R$ is a ring with sufficient conditions for equicodimensionality to hold. Since it has been shown afterwards that these hypotheses are superfluous for the main dimension formula,  so it seemed that all the corollaries would dispense as well with any extra hypotheses on $R$ (unless otherwise explicit in the statements).
\end{Remark}

Particular cases of interest in this work are as follows:

\begin{Corollary}
Let $I\subset R$ denote an ideal of grade $\geq 2$ in a Noetherian equicodimensional ring $R$ of dimension $3$. Let $I$ be generated by $r\geq 4$ elements. Then the following are equivalent:

{\rm (i)} $\mu(I_P)\leq r-1$ for every height $2$ prime $P\in \Spec R$

{\rm (ii)} $\dim \mathcal{S}_R(I)=\max \{4,r\}.$
\end{Corollary}
\demo (i) $\Rightarrow)$ (ii).
This implication follows from Proposition~\ref{Valla's bound} and does not require the bound $r\geq 4$.

(ii) $\Rightarrow)$ (i).
Let $\Ht P=2$. Then $\dim R/P=1$ since $R$ is equicodimensional.
Therefore,
$$1+\mu(I_P)=\dim R/P+\mu(I_P)\leq \dim \mathcal{S}_R(I)= \max \{4,r\}=r,$$
hence $\mu(I_P)\leq r-1$.
\qed

\begin{Example}\rm
The following simple example illustrates the above corollary.
Let $R=k[x,y,z]$ ($k$ a field) and let $I=(x,y)^3\cap (x,z)^3$, a $4$-generated ideal.
Then the number of generators of $I$ locally at the minimal primes $(x,y)$ and $(x,z)$ does not drop, and in fact, as it turns out, $\dim \mathcal{S}_R(I)=5> 4$.
\end{Example}

\begin{Corollary}
Let $R$ be a Noetherian equicodimensional, equidimensional catenary ring ane let $\phi$ denote an $n\times (n-1)$ matrix with entries in $R$. Suppose that $\dim R:=\mathfrak{d}<n$ and let $\mathcal{L}=I_1(\yy.\phi)$ denote a presentation ideal of $\mathcal{S}_R(I)$.
If $\phi$ satisfies condition $G_{\mathfrak{d}-1}$ and $I_1(\phi)$ has height $\mathfrak{d}$ then $\Ht \mathcal{L}=\mathfrak{d}$ {\rm (}maximal possible{\rm )}.
\end{Corollary}
\demo Since $n>\mathfrak{d}$, by Proposition~\ref{Valla's bound} one has $\dim \mathcal{S}_R(I)=n$.
By the assumptions on $R$, one has $\Ht \mathcal{L}=\mathfrak{d}+n-n=\mathfrak{d}$.
\qed

\medskip

If $I=(f_1,\ldots,f_r)\subset R$ is an equigenerated ideal in a standard graded ring over a field $k$ then, by a previous identification, the  kernel $Q$ of the surjection of graded $k$-algebras
$$k[y_1,\ldots,y_r]\to k[f_1t,\ldots,f_rt],\quad y_i\to f_it$$
is referred to as the presentation ideal of $\mathcal{F}_R(I)$. Obviously, $QR[y_1,\ldots,y_r]\subset \mathcal{J},$ where the latter is the presentation ideal of $\mathcal{R}_R(I)$  on $R[y_1,\ldots,y_n]$.
In the sequel we will omit the subscript $R$ in the notation of the special algebra to avoid confusion with the ground field $k$ over which it is naturally an algebra.

 As an immediate consequence of \cite[Theorem 5.1]{Tchernev} we have

\begin{Proposition}\label{Tchernev}
Let $R$ be a Noetherian Cohen-Macaulay ring and $I$ an ideal of $R$ such that ${\rm pd\,}I\leq 1.$ If $I$ satisfies $G_s$ then $\mathcal{S}_R(I)_i\simeq I^i$ for each $1\leq i\leq s-1.$	In particular, $\mathcal{A}_i=\{0\}$ for each $1\leq i\leq s-1.$
\end{Proposition}

A consequence of Proposition~\ref{Tchernev} is 

\begin{Corollary}\label{Tchernev_Cor} Let $R$ be a standard graded Cohen--Macaulay ring  over a field and let $I$ be an equigenerated perfect homogeneous ideal of $R$ of height $2$. If $I$ satisfies property $G_s$ then ${\rm indeg}(Q)\geq s.$ 
\end{Corollary}
Beyond the case where $s=\dim R$, which is a typical assumption in many known results, the knowledge of smaller cases can also be useful.
One situation is checking the birationality of a rational map via the main criterion of \cite{AHA}. Next is one example.

\begin{Example}\rm Let $I\subset R=k[x_1,x_2,x_3,x_4]$ be the height two perfect ideal generated by the maximal minors of the matrix
	$$\phi=\left[\begin{array}{ccccc}
	x_1&0&0&0&h_1\\
	x_2&x_1&0&0&h_2\\
	x_3&x_2&x_1&0&h_3\\
	0&x_3&x_2&0&h_4\\
	0&0&x_3&x_1&h_5\\
	0&0&0&x_4&x_3^{d}
	\end{array}\right],$$
	where $d\geq 1$ is a integer and $h_1,\ldots,h_5$ are forms of degree $d.$ 
	
Clearly, $I$ has one-less than maximal linear rank. Besides, without some calculation, one cannot decide whether $\ell(I)=4$, that is, if the map $\mathfrak{F}$ defined by the maximal ideals of $\phi$ is finite. Thus, \cite[Theorem 3.2]{AHA} is far from applicable to deduce that the map is birational onto the image.
Instead, we go about the birationality without deciding a priori if the map is finite.
Namely, a straightforward verification shows that a prime ideal $P$ containing $I_4(\phi)$ contains  $x_1^4,x_3^{d+3},x_2^3x_4$, hence, $\Ht I_4(\phi)\geq3.$  Thus, the ideal $I$ satisfies $G_3$. 

On the other hand, using the three rightmost linear syzygies of $\phi$ we  obtain a submatrix of the Jacobian dual matrix of $\mathfrak{F}$ as follows:
	$$\psi=\left[\begin{array}{ccccc}
	y_2&y_3&y_4&0\\
	y_3&y_4&y_5&0\\
	y_5&0&0&y_6
	\end{array}\right],$$
where the $y_i$'s are new variables over $k$. By Corollary~\ref{Tchernev_Cor}, the $3$-minor $y_5(y_3y_5-y_4^2)$  of $\psi$ does not vanish modulo $Q.$ Hence, $\psi$ has maximal rank over $Q$. By \cite[Theorem 2.18 (b)]{AHA}, $\mathfrak{F}$ is birational onto its image. Note, however, that  $I$ does not satisfy $G_4$ since $x_4$ appears at most on two columns, forcing $I_3(\phi)$ to have height $3$ only -- in addition, when $d\geq 2$, $G_4$ would imply that $\deg(\mathfrak{F})=d\geq 2$, the product of the syzygy degrees (see \cite[Corollary 3.2]{Cid-Ruiz}).
\end{Example}

\section{On the special fiber}

 \subsection{Relation to the associated graded ring}

As a piece of additional notation, let $I\subset R$ be a homogeneous ideal in a standard graded ring $R$ over a field and set $\fm:=(R_+)\subset R$. The image of an element $f\in I$ in $I/I^2\subset {\rm gr}_I(R)$ will be denoted by $f^o$.
In addition, set ${\rm gr}_I(R)_+:=(I/I^2){\rm gr}_I(R)$ for the ideal generated in positive degrees.

Similarly, the image  of an element $f\in I$ in $I/\fm I\subset \mathcal{F}(I)$) will be denoted by $f^{\ast}$. However, in this case, the notation will become obsolete throughout due to the identification of graded $k$-algebras $\mathcal{F}(I)\simeq k[I_{d}]=\bigoplus_{n\geq 0}[I^n]_{nd}\subset R$ established in (\ref{fiber_is_subalgebra}), whereby $f_i^{\ast}$ gets identified with $f_i\in [I]_{d}\subset  k[I_{d}]$.

We will use the following lemmata, certainly known in one way or another. Proofs are given for completeness.

\begin{Lemma}\label{characterization_CM}
	Let $R$ denote a standard graded ring  over an infinite field $k$ and let $I\subset R$ be a $d$-equigenerated homogeneous ideal. Let $J\subset I$ be a minimal reduction.
	One has:
	\begin{enumerate}
		\item[{\rm(i)}] $\mathcal{F}(I)\simeq k[I_{d}]$ is Cohen--Macaulay if and only if $k[I_{d}]$ is a free $k[J_d]$-module via the natural extension $k[J_d]\subset k[I_{d}]$.
		\item[{\rm(ii)}] If $\mathcal{F}(I)\simeq k[I_{d}]$ is Cohen--Macaulay then $\{1\}\cup B_1\cup\ldots\cup B_{\mathfrak{r}-1}$ is a free basis of $k[I_d]$ over $k[J_d]$, 
		where $B_n$ is the lifting of a homogeneous basis of the $k$-vector space $[I^n/JI^{n-1}]_{nd}$ and $\mathfrak{r}$ denotes the reduction number of $I$.
	\end{enumerate}
\end{Lemma}
\demo (i) This is the well-known argument, as applied to the  graded Noether normalization $k[J_d]\subset k[I_{d}]$.

(ii) Since $JI^{\mathfrak{r}}=I^{\mathfrak{r}+1}$, it is clear that the stated set generates $k[I_d]$ as $k[J_d]$-module.
By the Krull--Nakayama lemma, this set is a minimal set of generators, hence it must be a free basis.
\qed 

\begin{Lemma}\label{very_basic}
With ${\rm gr}_I(R)_+$ as above, one has:
\begin{enumerate}
	\item[{\rm (1)}] If $R$ is Cohen--Macaulay or a domain {\rm (}equidimensional catenary suffices{\rm )} then 
	$$\Ht {\rm gr}_I(R)_+\leq \Ht I.$$
	\item[{\rm (2)}] If ${\rm gr}_I(R)$ is Cohen--Macaulay then $\Ht {\rm gr}_I(R)_+\geq \Ht I$.
	\item[{\rm (3)}] If both $R$ and ${\rm gr}_I(R)$ are Cohen--Macaulay then $\Ht {\rm gr}_I(R)_+ = \Ht I$.
\end{enumerate}
\end{Lemma}
\demo (1) and (2) come out from the isomorphism ${\rm gr}_I(R)/{\rm gr}_I(R)_+\simeq R/I$, while (3) is the conjunction of (1) and (2).
\qed

\begin{Lemma}\label{regular_sequence_fiber}
	Let $R$ denote a standard graded ring  over an infinite field and let $I\subset R$ be a $d$-equigenerated homogeneous ideal of height $g$.  Let $f_1,\ldots,f_g\in I_{d}$ be a maximal regular sequence in $I$ such that 
	\begin{equation}\label{VV_condition}
(f_1,\ldots,f_g)\cap I^{n}=(f_1,\ldots,f_g)I^{n-1}.
	\end{equation}
	for every $n\geq 1.$ Then:
	\begin{enumerate}
		\item[{\rm (i)}] $f_1^{\ast},\ldots,f_g^{\ast}$ is a regular sequence in $\mathcal{F}(I);$  in particular, $\depth \mathcal{F}(I)\geq g.$
		\item[{\rm (ii)}] $(f_1,\ldots,f_g)_P$ is a minimal reduction of $I_P$ for every prime $P\supset I$ with $\Ht(P)=\Ht(I).$
		\item[{\rm(iii)}]  The sequence $f_1,\ldots,f_g$ can be extended to a sequence  generating a minimal reduction of $I.$ 
	\end{enumerate} 
\end{Lemma}
\demo (i) Having identified $\mathcal{F}(I)\simeq k[I_{d}]=\bigoplus_{n\geq 0}[I^n]_{nd}\subset R$, the assertion is merely a property of the graded inclusion $k[I_d]\subset R$, whereby one has to check that, for $0\leq i \leq g-1$, $f_{i+1}$ is a regular element on $\bigoplus_{n\geq 0} [I^n/(f_1,\ldots,f_i)I^{n-1}]_{nd}.$
Thus, let there be given $n$ and $h\in [I^{n}]_{nd}$ such that 
$hf_{i+1}\in (f_1,\ldots,f_i)[I^{n}]_{nd}\subset (f_1,\ldots,f_i)R$.

Since $\{f_1,\ldots,f_{i+1}\}$ is a regular sequence in $R$, then $h\in (f_1,\ldots,f_i)R$, hence  $h\in [(f_1,\ldots,f_i)R\cap I^{n}]_{nd}$.
Then, by \cite[Corollary 2.7]{VaVa}, $h\in [(f_1,\ldots,f_i)I^{n-1}]_{nd},$
as was to be shown.

 (ii) Localizing (\ref{VV_condition}) at a prime $P\supset I$ gives
$$(f_1,\ldots,f_g)_P\cap I_P^{n}=(f_1,\ldots,f_g)_PI_P^{n-1}.$$
If, moreover, $\Ht P=\Ht I$ then, since $\Ht (f_1,\ldots,f_g)=\Ht I$, it follows that $I_P$ is contained in the radical of $(f_1,\ldots,f_g)_P$.
Therefore, $I_P^n\subset (f_1,\ldots,f_g)_P$ for some $n>\!\!>0$, and hence, $I^{n}_P=(f_1,\ldots,f_g)_PI_P^{n-1}.$
This shows that  $(f_1,\ldots,f_g)_P$ is a reduction of $I_P$, hence a minimal one since any minimal reduction is minimally generated by $\ell(I_P)\geq \Ht I_P=g$ elements.

\smallskip

(iii) Since $k$ is infinite, and, by (i), $f_1,\ldots,f_g$  is a homogeneous regular sequence in $k[I_d]$, then it is a subset of a homogeneous system of parameters of $k[I_d]$. 
The latter generates a minimal reduction of the irrelevant maximal ideal of $k[I_d]$.
By a well-known fact, this reduction lifts to a minimal reduction of $I$ extending the $R$-regular sequence $f_1,\ldots,f_g$.
\qed

\smallskip

By \cite{VaVa}, one knows that the hypotheses of the previous lemma imply (actually, are equivalent to the assertion) that the residues $f_1^{o},\ldots,f_g^{o}$ constitute a regular sequence in  ${\rm gr}_I(R)$. The next lemma gives a condition on ${\rm gr}_I(R)$ for as to when such sequences exist at all.

\begin{Lemma}\label{general_sequence}
	Let $R$ denote a standard graded ring  over an infinite field and let $I\subset R$ be a $d$-equigenerated homogeneous ideal of height $g$.  If  the associated graded ring ${\rm gr}_I(R)$ is Cohen-Macaulay then there exist forms $f_1,\ldots,f_g\in I_{d}$ such that  $f_1^{o},\ldots,f_g^{o}$ is a regular sequence in  ${\rm gr}_I(R);$  if, in addition, $R$ is Cohen--Macaulay, then ${\rm gr}_I(R)_+$ has grade $g$.
\end{Lemma}
\demo It is basically a prime avoidance argument.
To wit,  letting $I=(\boldsymbol{\mathfrak{f}}_1,\ldots,\boldsymbol{\mathfrak{f}}_r)$, one has
\begin{eqnarray}
\Ht(\boldsymbol{\mathfrak{f}}_1^{o},\ldots,\boldsymbol{\mathfrak{f}}_r^{o})&=& \dim {\rm gr}_I(R)-\dim  \frac{{\rm gr}_I(R)}{(\boldsymbol{\mathfrak{f}}_1^{o},\ldots,\boldsymbol{\mathfrak{f}}_r^{o})}, \quad \mbox{since ${\rm gr}_I(R)$ is Cohen--Macaulay}\nonumber\\
&=&\dim R-\dim R/I\geq g.\nonumber
\end{eqnarray}
Now, since ${\rm gr}_I(R)$ is Cohen--Macaulay then $\Ht \mathfrak{P}=0$ for every associated prime $\mathfrak{P}$ of ${\rm gr}_I(R).$ Since $\Ht (\boldsymbol{\mathfrak{f}}_1^o,\ldots,\boldsymbol{\mathfrak{f}}_r^{o})\geq g> 0$  then 
$(\boldsymbol{\mathfrak{f}}_1^o,\ldots,\boldsymbol{\mathfrak{f}}_r^{o})$ is not contained in any of these primes.
By prime avoidance, a sufficiently general $k$-linear combination $f_1^0$ of $\boldsymbol{\mathfrak{f}}_1^o,\ldots,\boldsymbol{\mathfrak{f}}_r^o$ is such that $f_1^o$ is not contained in any  $\mathfrak{P}$, hence is ${\rm gr}_I(R)$-regular.
Let $f_1\in I_{d}$ be any lifting of $f_1^0$.

The general inductive is similar, whereby one assumes the existence of $f_1,\ldots, f_i\in I_{d}$, for $1\leq i\leq g-1$, such that  $f_1^o,\ldots,f_{i}^o$ is a regular sequence in ${\rm gr}_I(R)$. Since ${\rm gr}_I(R)/(f_1^o,\ldots,f_{i}^o)$ is Cohen-Macaulay then again, 
$$ \Ht(\boldsymbol{\mathfrak{f}}_1^{o},\ldots,\boldsymbol{\mathfrak{f}}_r^{o}){\rm gr}_I(R)/(f_1^o,\ldots,f_{i}^o)=g-1\geq 1,$$
hence, by prime avoidance, there exists a form $f_{i+1}\in I_{d}$ such that $f_{i+1}^{o}$ is regular on ${\rm gr}_I(R)/(f_1^o,\ldots,f_{i}^o)$.

The supplementary assertion is a consequence of Lemma~\ref{very_basic} (3).
\qed

\begin{Remark}\label{general_forms_OK}\rm
Although Lemma~\ref{general_sequence} claims the existence of {\em some} appropriate $d$-forms filling the request, it is clear from the proof that by choosing general forms in $I_{d}$ the result is equally obtainable.
\end{Remark}

The following example, which we don't claim to be a `smallest' or `simplest', shows that even general forms may not work if ${\rm gr}_I(R)$ is not Cohen--Macaulay.

\begin{Example}\rm
Let $I=(z^6,yz^5, y^2z^4, xy^2z^3, x^2y^2z^2, x^3y^3)\subset R=k[x,y,z]$, a height two perfect ideal of analytic deviation one.
A calculation with \cite{M2} gives that both gr$_I(R)$ and $\mathcal{F}(I)$ are almost Cohen--Macaulay (i.e., both have depth $2$).
If $f_1,f_2\in I_6$ are general forms, an additional computation outputs  $(f_1,f_2)\cap I^2\not\subset (f_1,f_2)I$ (alternatively, the grade of gr$_I(R)_+$ is $1$).  
\end{Example}

\begin{Corollary}\label{conjecture4equimultiple}
	Let $R$ denote a standard graded ring over an infinite field and let $I\subset R$ be a $d$-equigenerated homogeneous ideal. Suppose that $I$ is equimultiple {\rm(}that is, $\Ht I=\ell(I)${\rm).} If ${\rm gr}_I(R)$ is Cohen--Macaulay then so is the special fiber $\mathcal{F}(I)$.
\end{Corollary}

\begin{Remark}\rm
In the last corollary, if $R$ is Cohen--Macaulay, the assumption that ${\rm gr}_I(R)$ is Cohen--Macaulay can be weakened by just assuming that there exist $\Ht I-1$ forms  em $I_{d}$ whose images in ${\rm gr}_I(R)$ constitute a regular sequence.
\end{Remark}

The role of the Castelnuovo--Mumford regularity of a standard graded ring $A$ over a field, denoted ${\rm reg}(A)$, has been largely treated -- see \cite{EG} and \cite{Trung1998}, `classically' tailored, and \cite{KPU}, of a more recent vintage.

\begin{Lemma}
Let $A$ denote a standard graded domain over an algebraically closed field $k$ such that ${\rm reg}(A)\leq 1$. Then either $A$ is Cohen--Macaulay with minimal multiplicity or else $\depth A\leq \dim A-2$.
\end{Lemma}
\demo We may assume that ${\rm reg}(A)= 1$ as otherwise $A$ is a polynomial ring over $k$.
Then, writing $A=S/I$, where $S$ is a polynomial ring over $k$, $A$ has a $2$-linear resolution over $S$.
Suppose as if it were, that $\depth A=\dim A-1$. 
Then, by \cite[Theorem 1.5]{HerSri} one has 
$$e(A)=\frac{p!}{(p-1)!}-\frac{b_{p}}{(p-1)!}< p,$$
where $p$ denotes the homological dimension of $A$ over $S$ and $b_p$ is the last Betti number of the resolution.
On the other hand, $p=\edim A-\depth A$.
Therefore, $e(A)< \edim A-\dim A+1={\rm ecod} A+1$.
Since $A$ is a domain over an algebraically closed field, this contradicts a well-known result (\cite{EG}).
It follows that $A$ is Cohen--Macaulay and, necessarily, has minimal multiplicity.
\qed

\begin{Example}\rm
The following  well-known example shows that the above result does not hold in general if $A$ is not a domain: $I=(x_1,x_2)\cap (x_3,x_4)\subset R=k[x_1,x_2,x_3,x_4]$. Here $\depth R/I=1=\dim R/I -1$ and $\reg R/I=1$.
\end{Example}

Next is an application in our context.
Recall that the {\em analytic deviation} of an ideal $I$ (in a local or standard graded ring) is the difference $\ell(I)-\Ht I$.

\begin{Corollary}
	Let $R$ denote a standard graded  domain over an algebraically closed field and let $I\subset R$ be a $d$-equigenerated homogeneous ideal of analytic deviation one.
	If ${\rm gr}_I(R)$ is Cohen--Macaulay and has regularity at most one, then  $\mathcal{F}(I)$ is Cohen-Macaulay with minimal multiplicity.
\end{Corollary}
\demo By \cite{Ooishi} (see also \cite[Theorem 1.1]{Trung1998}),  $\reg {\rm gr}_I(R)=\reg \mathcal{R}(I)$, and by  \cite[Theorem 1.3 and Theorem 3.2]{RossiTrungs} or \cite[Section 1]{EU2012},  $\reg \mathcal{F}(I)\leq \reg\mathcal{R}(I).$ 
We may assume that $\reg \mathcal{F}(I)=1$ as otherwise there is nothing to prove. 
Since $\dim \mathcal{F}(I)=\ell(I)=\Ht I+1$ by the analytic deviation assumption, and $\depth \mathcal{F}(I)\geq \Ht I$, it follows that $\depth \mathcal{F}(I)\geq \dim  \mathcal{F}(I)-1$.
On the other hand, as $R$ is a domain and $\mathcal{F}(I)\simeq k[I_dt]\subset R[It]\subset R[t]$, then $\mathcal{F}(I)$ is a domain.
Therefore, the result is a consequence of the previous lemma.
\qed

\subsection{Analytically tight sequences of forms}

	Let $R$ denote a standard graded ring  over an infinite field and let $I\subset R$ be a  $d$-equigenerated homogeneous ideal. 
	We introduce the main notion of this section:

\begin{Definition}\label{super_general}\rm
	Let $I\subset R$ be a $d$-equigenerated homogeneous ideal in a standard graded ring over an infinite field.
	Set $\ell=\ell(I)$ (analytic spread) and $\mu=\mu(I)$ (minimal number of generators).
	We say that a sequence of forms $f_1,\ldots,f_{\ell}\in I_{d}$  is {\em analytically tight in power} $n$ if
	\begin{equation}\label{tight}
	[((f_1,\ldots,f_{\ell-1}):f_{\ell})\cap I^n]_{nd}=[(f_1,\ldots,f_{\ell-1})\cap I^n]_{nd}.
	\end{equation}
	The sequence $f_1,\ldots,f_{\ell}\in I_{d}$  is {\em analytically tight} if it is  analytically tight in power $n$ for every $n\geq 1.$ 
\end{Definition}
A few properties of this notion are:
\begin{enumerate}
	\item If $f_1,\ldots,f_{\ell}\in I_{d}$  is analytically tight in power $n$ then it is so in any power $\geq n$.
	\item If an equality as (\ref{tight}) holds for every subset $f_1,\ldots, f_i$, with $0\leq i\leq \ell$,  then we would be speaking of a {\em fully tight} sequence in power $n$. 
	\item A fully tight sequence satisfying in addition the Valabrega--Valla condition 
	$$(f_1,\ldots, f_i)\cap I^n\subset (f_1,\ldots, f_i)I^{n-1}$$ 
	is a {\em sequence of regular type} $n-1$ in the terminology of \cite{Trung1998}.
\end{enumerate}

The following result gives a measure of the relevance of the concept, at least in the case of ideals of analytic deviation one.

\begin{Theorem}\label{fiber_CM}
	Let $R$ denote a standard graded ring  over an infinite field and let $I\subset R$ be a  $d$-equigenerated homogeneous ideal of analytic deviation one and height $g$. Let  $f_1,\ldots,f_{g},f_{g+1}$ be forms in $I_{d}$ such that 
	\begin{enumerate}
		\item[{\rm(a)}] The images of $f_1,\ldots,f_{g}$ in ${\rm gr}_I(R)$ form a regular sequence.
		\item[{\rm(b)}] $J=(f_1,\ldots,f_{g+1})$ is a reduction of $I$.
	\end{enumerate}
	Then the following are equivalent:
	\begin{enumerate}
		\item[{\rm(i)}]  $f_1,\ldots,f_g,f_{g+1}$ analytically tight.
		\item[{\rm(ii)}] The special fiber $\mathcal{F}(I)$ is Cohen-Macaulay.
	\end{enumerate}
\end{Theorem} 
\demo 
(i)$\Rightarrow$ (ii) By Lemma~\ref{regular_sequence_fiber} (i), 
$f_1,\ldots,f_{g}$ is a regular sequence in 
$$\mathcal{F}(I)=k[I_d]=\bigoplus_{n\geq 0} [I^n]_{nd},$$ regarded as elements of $[I]_d$.
Thus,  we are to prove that $f_{g+1}\in [I]_{d}$ is a non-zero-divisor of
$\bigoplus_{n\geq 0} [I^n/(f_1,\ldots,f_i)I^{n-1}]_{nd}$.

So, let $h\in [I^{n_0}]_{n_0d}$ be such that $hf_{g+1}=0$ in $\bigoplus_{n\geq 0} [I^n/(f_1,\ldots,f_i)I^{n-1}]_{nd}$. This means that
$h\in [((f_1,\ldots,f_g)I^{n_0}:f_{g+1})\cap I^{{n_0}}]_{{n_0}d}\subset [((f_1,\ldots,f_g):f_{g+1})\cap I^{{n_0}}]_{{n_0}d}.$
Therefore,
$$\begin{array}{ccclc}
h&\in & [(f_1,\ldots,f_g)\cap I^{n_0}]_{{n_0}d}& \mbox{because $f_1,\ldots,f_g,f_{g+1}$ is analytically tight }\\
&=&[(f_1,\ldots,f_g)I^{{n_0}-1}]_{{n_0}d}& \mbox{because $f_1^o,\ldots,f_g^o$ is a regular sequence in gr$_I(R)$}. 
\end{array}.$$

(ii)$\Rightarrow$(i) Since $J=(f_1,\ldots,f_{g+1})$ is a reduction of $I$ then $f_1,\ldots,f_{g+1}$, regarded as elements of $[I]_d$, is a system of parameters of $\mathcal{F}(I)=k[I_d]$, hence a regular sequence of $\mathcal{F}(I)$ since the latter is assumed to be Cohen-Macaulay.

The proof is similar to the argument in the previous item, only changing sense.
Thus, for the essential inclusion to be seen, pick $h\in [((f_1,\ldots,f_{g}):f_{g+1})\cap I^n]_{nd}.$ Then 
$$hf_{g+1}\in [(f_1,\ldots,f_{g})\cap I^{n+1}]_{(n+1)d}=
[(f_1,\ldots,f_{g}) I^{n}]_{(n+1)d},$$
where the equality results from assumption (a) via \cite[Corollary 2.7]{VaVa}.
 Since $f_{g+1}$ is a non-zero-divisor of
 $\bigoplus_{n\geq 0} [I^n/(f_1,\ldots,f_i)I^{n-1}]_{nd}$, one has
\begin{eqnarray*} 
		h&\in&[(f_1,\ldots,f_{g}) I^{n-1}]_{nd}\nonumber\\
	&=&[(f_1,\ldots,f_g)\cap I^n]_{nd} \quad \mbox{because $f_1^o,\ldots,f_g^o$ is a regular sequence in gr$_I(R)$}\nonumber.\quad\mbox{\qed}
\end{eqnarray*}

In the next proposition sufficient conditions are established for a sequence to be analytically tight in power $n$ for every $ n >\!\!>0$.
The proof is an adaptation of an argument in \cite[Remark 2.1]{HucHun}.

\begin{Proposition}\label{asymptoticallyissuper}
		Let $R$ denote a standard graded ring  over an infinite  field $k$ and let $I\subset R$ be an unmixed $d$-equigenerated homogeneous ideal of analytic deviation one and height $g$. Suppose that $\{f_1,\ldots,f_{g},f_{g+1}\}\subset [I]_d$ generates a reduction of $I$, with $f_1,\ldots,f_g$ a regular sequence in $R$.
	One has:
		\begin{enumerate}
			\item[{\rm(i)}] $f_1,\ldots,f_g,f_{g+1}$ is analytically tight in sufficiently high powers.
			\item[{\rm(ii)}] If, moreover,  $(f_1,\ldots,f_g)\cap I^n=(f_1,\ldots,f_g)I^{n-1}$, for every $n\geq 1$. Then  $f_1,\ldots,f_g,f_{g+1}$ is analytically tight in powers  $n\geq n_0,$ where $n_0=\max\{r(I_P)\,|\,P\in{\rm Min}(R/I)\}+1.$ 
		In particular, in the case where  $I$ is generically a complete intersection, $f_1,\ldots,f_g,f_{g+1}$ is analytically tight in all powers.
	\end{enumerate}
\end{Proposition}
\demo (i)  By degree reasons, it suffices to show that  there is an integer $n_0$ such that 
$$((f_1,\ldots,f_g):f_{g+1})\cap I^n=(f_1,\ldots,f_g)\cap I^{n},$$
for every $n\geq n_0.$ 

 To prove the essential inclusion set $\ass(R/(f_1,\ldots,f_g))=\{P_1,\ldots,P_r,Q_1,\ldots,Q_s\}$, where $\ass(R/I)=\{P_1,\ldots,P_r\}.$
Since $I$ is assumed to be unmixed, $\Ht P_i=\Ht Q_j=g$ for all $i,j$. Write $(f_1,\ldots,f_g)=N\cap N_1$, where $N$ is the intersection of the primary components of $(f_1,\ldots,f_g)$ corresponding to $\{P_1,\ldots,P_s\}$ and $N_1$ is the intersection of the primary components of $(f_1,\ldots,f_g)$ corresponding to $\{Q_1,\ldots,Q_s\}.$ Since $\sqrt{I}=\sqrt{N}$,  there exists an integer $n_0$ such that $I^n\subset N$ for every $n\geq n_0.$ 

Pick $h\in ((f_1,\ldots,f_g):f_{g+1})\cap I^n$ with $n\geq n_0.$ Then $hf_{g+1}\in N\cap N_1$ and $h\in I^{n}.$ Since $\sqrt{(f_1,\ldots,f_g,f_{g+1})}=\sqrt{I}$ and $I\not\subset Q_j$ for $j=1,\ldots,s,$ forcefully $f_{g+1}\notin Q_{j}$ for $j=1,\ldots,s.$ Thus, $h\in N_1\cap I^n\subset N_1\cap N=(f_1,\ldots,f_g)$, as required.

(ii)  As in the proof of item (i), set $\ass(R/(f_1,\ldots,f_g))=\{P_1,\ldots,P_r,Q_1,\ldots,Q_s\}$, where $\ass(R/I)=\{P_1,\ldots,P_r\}$ and $(f_1,\ldots,f_g)=N\cap N_1$, where $N$ is the intersection of the primary components of $(f_1,\ldots,f_g)$ corresponding to $\{P_1,\ldots,P_s\}$ and $N_1$ is the intersection of the primary components of $(f_1,\ldots,f_g)$ corresponding to $\{Q_1,\ldots,Q_s\}.$ By a similar token as in the proof of (i), it is sufficient to prove that $I^{n_0}\subset N$ for $n_0=\min\{r(I_P)\,|\,P\in{\rm Min}(R/I)\}+1$. 

Now, by assumption, $(f_1,\ldots,f_g)\cap I^n=(f_1,\ldots,f_g)I^{n-1}$, for every $n\geq 1.$ In particular, for every $P\in \ass(R/I)={\rm Min}(R/I),$ the images of $f_1,\ldots,f_g$ in $R_P$ is a regular sequence and $$(f_1,\ldots,f_g)_P\cap I_P^n=(f_1,\ldots,f_g)_PI_P^{n-1}$$
for every $n\geq 1.$ Thus, the images of $f_1,\ldots,f_g$ in ${\rm gr}_{I_P}(R_P)$ form a regular sequence. In particular, $\depth {\rm gr}_{I_P}(R_P)=\dim {\rm gr}_{I_P}(R_P)=g,$ that is, ${\rm gr}_{I_P}(R_P)$ is Cohen-Macaulay.  Hence, by \cite[Theorem 1.2]{Trung1987}, the reduction number $r(I_P)$ of $I_P$ is independent of the minimal reduction. 
By  Lemma~\ref{regular_sequence_fiber} (ii), $(f_1,\ldots,f_g)_P$ is  a minimal reduction of $I_P$, hence $r(I_P)\leq n_0-1$ by assumption. Consequently, $$I_P^{n_0}=(f_1,\ldots,f_g)_PI^{n_0-1}_P\subset (f_1,\ldots,f_g)_P=N_{P}.$$
But, since $I_P^{n_0}\subset N_P$ for every $P\in \ass(R/N)=\ass(R/I)$, we have  $I^{n_0}\subset N$ as desired.

To see the supplementary assertion, if  $I$ is generically a complete intersection then $ I_P = (f_1, \ldots, f_g) _P$ over any minimal prime $P$ of $I$. Hence, in this case, the argument shows that one can take $n_0=1.$
\qed

\begin{Remark}\rm
Recall again that, by \cite[Corollary 2.7]{VaVa}, the additional hypothesis of item (ii), along with the standing assumption of item (a), makes it  equivalent to requiring that the images of $f_1,\ldots,f_g$ in ${\rm gr}_I(R)$ are a regular sequence. Computationally, this angle is preferable.
\end{Remark}

\begin{Corollary}
Let $R$ denote a standard graded Cohen--Macaulay ring  over an infinite field and let $I\subset R$ be an equigenerated homogeneous ideal satisfying the following conditions:
\begin{enumerate}
	\item[{\rm (i)}] $I$ is unmixed and generically a complete intersection.
	\item[{\rm (ii)}] $I$ has analytic deviation one.
	\item[{\rm (iii)}] ${\rm gr}_I(R)$ is Cohen--Macaulay.
\end{enumerate}

 Then $\mathcal{F}(I)$ is Cohen--Macaulay.
\end{Corollary}

An example of the above corollary in dimension three is given by the ideal generated by the $(m-1)$-products of $m\geq 3$ linear forms in $R$ defining a central generic arrangement
This is a special case of the main result in \cite{ProudSpei} (see also \cite[Proposition 4.1]{GST}).

\begin{Proposition}\label{getting_anal_tight}
	Let $R$ denote a standard graded  ring  over an infinite field  and let $I\subset R$ be a $d$-equigenerated homogeneous ideal of analytic deviation one and of height $g$. 
	\begin{enumerate}
		\item[{\rm(i)}] If ${\rm gr}_I(R)$ is Cohen--Macaulay, there exists  a sequence of forms $f_1,\ldots,f_g,f_{g+1}$ in $I_{d}$ such that $f_1^{o},\ldots,f_g^{o}$ is a regular sequence in ${\rm gr}_I(R)$ and $J=(f_1,\ldots,f_{g+1})$ is a minimal reduction of $I.$
		\item[{\rm (ii)}] Suppose, moreover, that $I$ is unmixed. If $R$ and  $\mathcal{R}_R(I)$ are Cohen--Macaulay locally at the minimal primes of $I$, then any sequence as in {\rm (i)}  is analytically tight in powers $\geq \ell(I)-1$.
		In particular, this is the case when $R$ is regular locally at the minimal primes of $R$ and ${\rm gr}_I(R)$ is Cohen--Macaulay locally at these primes.
	\end{enumerate}
\end{Proposition}
\demo (i) This is an immediate consequence of Lemma~\ref{general_sequence} and Lemma~\ref{regular_sequence_fiber}.

(ii) Since both $R_P$ and $\mathcal{R}_{R_P}(I_P)$ are Cohen--Macaulay for every $P\in {\rm Min}(R/I)$, then \cite[Theorem 5.1]{AHT1995} yields $r(I_P)\leq \dim R_P-1=g-1$ for every such $P$. Hence, in the notation of Proposition~\ref{asymptoticallyissuper}, 
$$n_0=\max\{r(I_P)\,|\,P\in {\rm Min}(R/I)\}+1\leq g=\ell(I)-1.$$
On the other hand, since $f_1^{o},\ldots,f_g^{o}$ is a regular sequence in ${\rm gr}_I(R)$, then it satisfies the Valabrega--Valla condition
$(f_1,\ldots,f_g)\cap I^n=(f_1,\ldots,f_g)I^{n-1}$, for every $n\geq 1$.
Since $I$ is assumed to be unmixed, Proposition~\ref{asymptoticallyissuper} (ii) says that the sequence $f_1,\ldots,f_g,f_{g+1}$ is analytically tight in powers $\geq \ell(I)-1$.

The supplementary statement is a consequence of the Brian\c con--Skoda theory (\cite[Proposition 7.3.35]{SiBook}).
\qed

\medskip

The existence of analytically tight sequences whose forms are general might be a fairly encompassing phenomenon, regardless of Cohen--Macaulay hypothesis on either $\mathcal{F}(I)$ or ${\rm gr}_I(R)$.
Computational evidence suggests that perfect height two ideals in dimension three may be an adequate source. Next is an example to keep this in mind.

\begin{Example}\rm Let $I=(z^6, yz^5, xyz^4, xy^2z^3, xy^3z^2, x^2y^3z, x^3y^3)\subset R=k[x,y,z]$, a perfect height two ideal with $\ell(I)=3$. For general $f_1,f_2,f_3\in I_{6},$ a calculation with \cite{M2} gives
	$$[((f_1,f_{2}):f_{3})\cap I^n]_{6n}=[(f_1,f_{2})\cap I^n]_{6n},$$
	for $n=1$, hence for every $n\geq 1$ as well.
	Hence, $\{f_1,f_2,f_3\}$ is analytically tight. However, although $r_{(f_1,f_2,f_3)}(I)={\rm reg}(\mathcal{F}(I))=2$, the special fiber $\mathcal{F}(I)$ is not Cohen-Macaulay. In confront with Theorem~\ref{fiber_CM}, one derives that the associated graded ring is not Cohen--Macaulay -- actually, $\depth {\rm gr}_I(R)=\depth {\rm gr}_I(R)_+=1$.
\end{Example}

\subsection{Analytically adjusted sets of forms}

Some of the ideas of the previous subsection can be stated in terms of certain other expected equalities.

\begin{Definition}\label{J-equality}\rm
	Let $I\subset R$ be a $d$-equigenerated homogeneous ideal in a standard graded ring over an infinite field.
	Set  $\mu=\mu(I)$ (minimal number of generators).
	A set $\{f_1,\ldots,f_l\}\subset I_{d}$ of $k$-linearly independent forms  will be said to be $l$-{\em adjusted}  if  $\mu(IJ)=l\mu(I)-{l \choose 2}$, where $J=(f_{1},\ldots,f_l)$.
	
	If $l=\ell(I)$ is the analytic spread of $I$, an $l$-adjusted set will also be said to be {\em analytically adjusted}.
\end{Definition}

\begin{Remark}\rm (1) Clearly, for any choice of such indices, one always has the inequality $\mu(IJ)\leq l\mu(I)-{l \choose 2}$.
	We observe that a natural candidate for an analytically adjusted set of generators is one such that $J\subset I$ is a minimal reduction --
	recall that, in the equigenerated case, one can take a reduction $J\subset I$ that is a subset of a minimal set of generators of $I$. 
	
	(2) An advantage of equigeneration is that one essentially counts vector space dimensions, hence $\mu(I^2)\geq \mu(IJ)$ in the setup of Definition~\ref{J-equality}.
	Therefore, the present notion is somewhat stronger than the Freiman kind of lower bound inequality (see \cite[Theorem 1.9]{HerHibiZhu}).
	The same paper gives examples where any reasonable expectation may go wrong if $I$ is not equigenerated, even in the monomial case.
\end{Remark}	

Next are two  situations in which the above condition is satisfied.

\begin{Proposition}\label{no_quadratics_is_adjusted}
	Let $I\subset R$ be a $d$-equigenerated homogeneous ideal in a standard polynomial ring over an infinite field. Suppose that the special fiber $\mathcal{F}(I)$ admits no quadratic polynomial relations -- e.g., height two perfect ideals satisfying $G_3$ and, more generally, syzygetic equigenerated ideals. Then, any set of $k$-linearly independent forms $f_1,\ldots,f_{l}\subset I_{d}$ is $l$-adjusted.
\end{Proposition}
\demo
Since $f_1,\ldots,f_l\in I_{d}$ are  $k$-linearly independent, one can extend them to a minimal set $\{f_1,\ldots,f_{\mu}\}$ of generators of $I$. The ideal $(f_1,\ldots,f_l)I$ is a $2d$-equigenerated homogeneous ideal generated by the following set of forms of degree $2d:$
$$B=\{f_1f_j\,|\,1\leq j\leq \mu(I)\}\cup\{f_2f_j\,|\,2\leq j\leq \mu(I)\}\cup\ldots\cup\{f_lf_j\,|\,i\leq j\leq \mu(I)\}$$
It suffices to prove that this set is $k$-linearly independent. 
Indeed, if so then the set $B$ is a basis of the $k$-vector space $[(f_1,\ldots,f_l)I]_{2d}.$ In particular, $B$ is a minimal set of homogneous generators of $(f_1,\ldots,f_l)I$. Since $B$ has
$\mu((f_1,\ldots,f_l)I)=l\mu(I)-{l \choose 2}$ elements, we are through.

But  a nonzero $k$-linear combination of the elements of $B$  yields a nonzero quadratic polynomial relation of the special fiber $\mathcal{F}(I)$, which contradicts the assumption. 
\qed
 
 \begin{Proposition}\label{regular_sequence_is_adjusted}
Let $R$ denote a standard graded ring over an infinite field $k$ and let $I\subset R$ be a  $d$-equigenerated homogeneous ideal.
If $f_1,\ldots,f_l\in [I]_d$ are forms such that $f_1,\ldots,f_l$ is a regular sequence in  $k[I_d]$ then, for any $1\leq i\leq l$, the subset $\{f_1,\ldots,f_i\}$ is  $i$-adjusted.
In particular,  if the special fiber $\mathcal{F}(I)$ is Cohen-Macaulay then every set of forms $\{f_1,\ldots,f_{\ell}\}\subset I_d$ generating a minimal reduction $J$ of $I$ is analytically adjusted.
 \end{Proposition}
\demo First extend the regular sequence $f_1,\ldots,f_l$  to a full $k$-vector basis $$\{f_1,\ldots,f_l,f_{l+1},\ldots,f_{\mu(I)}\}$$
of $I_{d}.$ Clearly, $\mu((f_1,\ldots,f_i)I) =\dim_k[(f_1,\ldots,f_i)I]_{2d}$.
As already pointed out earlier, the inequality $\dim_k[(f_1,\ldots,f_i)I]_{2d}\leq i\mu(I)-{i\choose 2} $ is immediate because $[(f_1,\ldots,f_i)I]_{2d}$ has an obvious set of generators with $i\mu(I)-{i\choose 2}$ elements. 
To prove the reverse inequality we induct on $i.$
For $i=1$ the result is clear. Now, for $1< i\leq l$ let $\mathfrak{B}$ denote a $k$-vector base of $[(f_1,\ldots,f_{i-1})I]_{2d}.$ A vanishing $k$-linear combination of elements in $\mathfrak{B}\cup \{f_i^2,f_{i}f_{i+1},\ldots,f_{i}f_{\mu(I)}\} $ can be written as
\begin{eqnarray}
f_{i}(\lambda_if_i+\cdots+\lambda_{\mu(I)}f_{\mu(I)})&=& \mbox{$k$-linear combination of the elements of }\mathfrak{B}\nonumber\\
&\in& (f_1,\ldots,f_{i-1})I\nonumber
\end{eqnarray}
Thus, since $f_i$ is $\mathcal{F}(I)/(f_1,\ldots,f_{i-1})$-regular we have $\lambda_if_i+\cdots+\lambda_{\mu(I)}f_{\mu(I)}=0$, hence $\lambda_i=\cdots=\lambda_{\mu(I)}=0$ because
$f_i,\ldots,f_{\mu(I)}$ are $k$-linearly independent.   Consequently, one has that $\mathfrak{B}\cup \{f_i^2,f_{i}f_{i+1},\ldots,f_{i}f_{\mu(I)}\} $ is a $k$-linearly independent subset of $[(f_1,\ldots,f_i)I]_{2d}$. Since its number of elements is $\dim_k[(f_1,\ldots,f_{i-1})I]_{2d}+\mu(I)-i+1 $, by induction,  $\dim_k[(f_1,\ldots,f_{i})I]_{2d}\geq i\mu(I)-{i\choose 2} $ as desired.

If $\mathcal{F}(I)$ is moreover Cohen--Macaulay, then any minimal set of generators of a minimal reduction is a regular sequence.
Still, an independent proof can be given as follows.
Clearly, the $k$-vector space $[JI]_{2d}$ is spanned by the set $\mathfrak{B}:=\{f_if_j|1\leq i<j\leq \ell\}\cup \{f_1h_i,\ldots,f_{\ell}h_i|1\leq i\leq  s\}$ of cardinality $\ell\mu(I)-{\ell\choose2}$.

One claims that $\mathfrak{B}$ is $k$-linearly independent.
Indeed, let 
\begin{equation}\label{linear_combination}
\sum_{1\leq i<j\leq \ell} a_{i,j} f_if_j 
+ \sum_{\tiny \begin{array}{l}
	1\leq j\leq \ell\\
	1\leq t\leq s
	\end{array}} b_{j,t}f_jh_i=0
\end{equation}
be a relation of $k$-linear dependence thereof.
Now,  $a_{i,j}f_if_j\in [J]_{2d}$ and $b_{j,t}f_j\in [J]_d$, for all $i,j,t$.
Since $\mathcal{F}(I)=k[I_d]$ is Cohen--Macaulay,  by Lemma~\ref{characterization_CM}, $k[I_d]$ is a free module over $k[J_d]$, with $\{1,h_1,\ldots,h_s\}$  a subset of a free basis of $k[I_d]$.
Therefore, \eqref{linear_combination} implies that $\sum_{1\leq i<j\leq \ell} a_{i,j} f_if_j=\sum_{1\leq j\leq\ell} b_{j,t}f_j=0$ for all choices of $i,j,t$. 
But, $\{f_1,\ldots,f_{\ell}\}$ is an algebraically independent set over $k$, thereby leading to $a_{i,j}=b_{j,t}=0$, for all $i,j,t$.

Therefore, $\mathfrak{B}$ is a basis for the $k$-vector space $[JI]_{2d},$ as was to be shown.
\qed

\smallskip

One has a  weak converse:

\begin{Theorem}\label{weak_converse}
	Let $R$ denote a standard graded ring over an infinite field  and let $I\subset R$ be a $d$-equigenerated homogeneous ideal such that
	\begin{enumerate}
		\item[{\rm(a)}] $I$ is unmixed of analytic deviation one.
		\item[{\rm(b)}] For any minimal prime $P$ of $R/I$ the reduction number of $I_P$ is at most one.
	\end{enumerate}
	If $f_1,f_2,\ldots,f_{g+1}$ is a sequence of forms in $I_{d}$ satisfying the properties in {\rm Proposition~\ref{getting_anal_tight}\, (i)}, then the following are equivalent:
	\begin{enumerate}
		\item[{\rm(i)}]  $\mathcal{F}(I)$ is Cohen-Macaulay.
		\item[{\rm(ii)}] $f_1,f_2,\ldots,f_{g+1}$ is analytically adjusted.
		\item[{\rm(iii)}] $f_1,f_2,\ldots,f_{g+1}$ is analytically tight.
	\end{enumerate}
\end{Theorem}
\demo 
The equivalence (i)$\Leftrightarrow$(iii)  is a consequence of Theorem~\ref{fiber_CM}, while the implication (i)$\Rightarrow$(ii) is a consequence of Proposition~\ref{regular_sequence_is_adjusted}. Therefore, it remains to prove that (ii)$\Rightarrow$(i).

By Lemma~\ref{regular_sequence_fiber} (i),  $\{f_1,\ldots,f_g\}\subset [I]_{d}$ is a regular sequence in $\bigoplus_{n\geq 0}[I^n]_{nd}$.
We will be through if $f_{g+1}\in [I]_{d}$ is regular on 
$\bigoplus_{n\geq 0} [I^n/(f_1,\ldots,f_g)I^{n-1}]_{nd}$.
For this, let $h\in [I^n]_{nd}$ be such that $hf_{g+1}\in [(f_1,\ldots,f_{g})I^{n}]_{(n+1)d}.$ 
We separate the argument in two cases:

\medskip

\noindent{\bf Case  $n=1$:} We can extend the set $\{f_1,\ldots,f_{g+1}\}$ to a basis $\{f_1,\ldots,f_{\mu(I)}\}$ of the $k$-vector space $I_{d}.$ Since $\{f_1,\ldots,f_{g+1}\}$ is analytically adjusted then the set
$$\mathfrak{B}=\{f_if_j|1\leq i\leq j\leq g+1\}\cup \{f_if_j|1\leq i\leq g+1,g+2\leq j\leq \mu(I)\}$$
is a basis of the $k$-vector space $[JI]_{2d}$, where $J=(f_1,\ldots,f_{g+1}).$ On the other hand, since $\{f_1,\ldots,f_{g}\}$ is a regular sequence in  $\mathcal{F}(I)=k[I_d]$, it is $g$-adjusted by Proposition~\ref{regular_sequence_is_adjusted}, i.e., $\mu((f_1,\ldots,f_{g+1})I)=g\mu(I)-{g \choose 2}$, and hence, the set 
$$\mathfrak{B}'=\{f_if_j|1\leq i\leq j\leq g\}\cup \{f_if_j|1\leq i\leq g,g+1\leq j\leq \mu(I)\}\subset \mathfrak{B}$$
is a basis of the $k$-vector space $[(f_1,\ldots,f_{g})I]_{2d}.$

Since $n=1$, we have $hf_{g+1}\in[(f_1,\ldots,f_{g})I]_{2d}$ by hypothesis.

Writing $h=\lambda_1f_1+\cdots+\lambda_{g+1}f_{g+1}+\cdots+\lambda_{\mu(I)}f_{\mu(I)}$, then
$$hf_{g+1}=\lambda_1f_1f_{g+1}+\cdots+\lambda_{g+1}f_{g+1}^2+\cdots+\lambda_{\mu(I)}f_{\mu(I)}f_{g+1}\in [(f_1,\ldots,f_{g})I]_{2d}$$
implies that $\lambda_{g+1}f_{g+1}^2+\cdots+\lambda_{\mu(I)}f_{\mu(I)}f_{g+1}$
is a $k$-linear combination  of the elements of $\mathfrak{B}'.$ But, $\mathfrak{B}=\mathfrak{B}'\cup \{f_{g+1}^2,\ldots,f_{\mu(I)}f_{g+1}\}$ and $\mathfrak{B}'\cap \{f_{g+1}^2,\ldots,f_{\mu(I)}f_{g+1}\}=\emptyset.$ Therefore, it must be the case that $\lambda_{g+1}=\cdots=\lambda_{\mu(I)}=0.$  
Hence, $h\in[(f_1,\ldots,f_{g})]_{d}.$

\medskip

\noindent{\bf Case  $n\geq2$:} 
By assumption, 
$hf_{g+1}\in [(f_1,\ldots,f_{g})I^{n}]_{(n+1)d}$, i.e., $h\in[((f_1,\ldots,f_g):f_{g+1})\cap I^{n}]_{nd}.$
By Proposition~\ref{asymptoticallyissuper} (ii), assumption (b) implies that $f_1,\ldots,f_{g+1}$ is an analytically tight sequence in powers  $\geq 2$, hence
$$h\in [(f_1,\ldots,f_g)\cap I^n]_{nd}=[(f_1,\ldots,f_g)I^{n-1}]_{nd},$$
where the equality is the Valabrega--Valla condition, as a consequence of $f_1^o,\ldots,f_g^o$ being a regular sequence in gr$_I(R)$.
\qed

\section{Applications in dimension three}

The following is a consequence of the last part of the previous section:

\begin{Corollary}\label{weak_converse_dim3}
	Let $R$ denote a standard graded ring of dimension three over an infinite field  and let $I\subset R$ be a perfect $d$-equigenerated homogeneous ideal of height $2$ and analytic spread $3$. 
	Suppose that, for any height two prime $P$,  the reduction number of $I_P$ is at most one.
	For any sequence of forms $f_1,f_2,f_{3}$ in $I_{d}$ satisfying the properties in {\rm Proposition~\ref{getting_anal_tight}\, (i)},  the following are equivalent:
	\begin{enumerate}
		\item[{\rm(i)}]  $\mathcal{F}(I)$ is Cohen-Macaulay.
		\item[{\rm(ii)}] $f_1,f_2,f_{3}$ is analytically adjusted.
		\item[{\rm(iii)}] $f_1,f_2,f_{3}$ is analytically tight.
	\end{enumerate}
\end{Corollary}
\demo 
It is a reading of Theorem~\ref{weak_converse} in dimension three.
\qed

\smallskip

\begin{Remark}\rm
We note that the assumption on the local reduction number is satisfied when both $R$ and $\mathcal{R}(I)$ are Cohen--Macaulay locally at the minimal primes of $R/I$.
\end{Remark}
For height two equigenerated ideals, not necessarily unmixed, one can file the following result:

\begin{Theorem}\label{degreeRees_CM}
	Let $R$ denote a $3$-dimensional standard graded Cohen--Macaulay domain over an infinite field and let $I\subset R$ be a perfect $d$-equigenerated homogeneous ideal of height $2$. Suppose that:
	\begin{itemize}
		\item[{\rm(i)}] $\mu(I)\geq 4$ and $\ell(I)=3.$
		\item[{\rm(ii)}] The Rees algebra of $I$ is Cohen--Macaulay.
		\item[{\rm(iii)}] The defining ideal of $\mathcal{F}(I)$ is generated in degrees $\geq 3$.
	\end{itemize} 
	Then:
	\begin{enumerate}
		\item[{\rm(a)}] $\mathcal{F}(I)$ is Cohen--Macaulay and the reduction number of $I$  is $2.$ 
		\item[{\rm(b)}]The defining ideal of $\mathcal{F}(I)$  is equigenerated in degree $3$ and the minimal graded free resolution of $\mathcal{F}(I)$  is $3$-linear.
		\item[{\rm(c)}] The  multiplicity  of $\mathcal{F}(I)$ is ${m-1\choose 2},$ where $m:=\mu(I)=\dim_k(I_{d})$. 
	\end{enumerate}
\end{Theorem}
\demo (a) The assumption in (ii) implies that the associated graded ring gr$_I(R)$ is Cohen--Macaulay. Then, by Lemma~\ref{general_sequence} and Lemma~\ref{regular_sequence_fiber} (iii), there exist forms $f_1,f_2,f_3$ in $I_{d}$ such that $f_1^o,f_2^o$ is a regular sequence in gr$_I(R)$ and $J=(f_1,f_2,f_{3})$ is a minimal reduction of $I.$  By Proposition~\ref{no_quadratics_is_adjusted}, condition (iii) implies that $\{f_1,f_2,f_3\}$   is an analytically adjusted set. Then the implication (ii) $\Rightarrow$ (i) of Theorem~\ref{weak_converse} (b) gives that $\mathcal{F}(I)$ is Cohen--Macaulay as well.

Since $\mathcal{F}(I)$ is Cohen--Macaulay,  the reduction number $r(I)$ is independent on the minimal reduction and coincides with the Castelnuovo-Mumford regularity $\reg(\mathcal{F}(I))$ (see, e.g., \cite[Proposition 1.2]{GST}).
On the other hand, since $I$ is minimally generated by at least four elements, condition (iii) implies that $\reg(\mathcal{F}(I))\geq 2.$ 
But by (ii),  $r(I)\leq \ell(I)-1=2$ (see \cite[Theorem 5.6]{AHT1995}).

(b) Since  $\reg(\mathcal{F}(I))=2$, then  (iii) implies that the defining ideal of $\mathcal{F}(I)$ is equigenerated in degree $3$. In particular, the minimal graded free resolution of $\mathcal{F}(I)$  is $3$-linear.

(c) Since the minimal graded free resolution of $\mathcal{F}(I)$ is $3$-linear, this follows immediately from \cite[Theorem 1.2]{HM}.
 \qed
 
\begin{Remark}\rm
Suitable versions of Theorem~\ref{degreeRees_CM}  for arbitrary number of variables may be stated. Both results may profit from the still standing Simis--Vasconcelos question as to whether, for an equigenerated homogeneous ideal $I$ in a standard polynomial ring over a field, if the Rees algebra of $I$ is Cohen--Macaulay then so is the special fiber of $I$ (see Corollary~\ref{conjecture4equimultiple} for the equimultiple case, valid under weaker hypothesis).
\end{Remark}

The following calculation is certainly well-known. A proof is given for the reader's convenience.

\begin{Lemma}\label{multiplicity_of_perfect_cod2}
Let $I\subset R=k[x,y,z]$ be a $d$-equigenerated homogeneous perfect ideal of height $2$ with minimal graded free resolution
$$0\to \bigoplus_{i=1}^{N-1} R(-d-m_i)\stackrel{\phi}\lar R(-d)^N\to R\to R/I\to 0.$$
Then $$e(R/I)=\frac{d^2+ \sum_{i=1}^{N-1}m_i^2}{2}.$$
\end{Lemma}
\demo
The Hilbert series of $R/I$ calculated from the given resolution is
$$H_{R/I}(t)=\frac{B_{R/I}(t)}{(1-t)^3}$$
where $$B_{R/I}(t)=1-Nt^{d}+\sum_{i=1}^{N-1}t^{d+m_i}.$$
By a well-known formula (see, e.g., \cite[Corollary 7.4.12]{SiBook})
\begin{eqnarray}
e(R/I)&=&\frac{1}{2}\,\frac{\partial^2B_{R/I}(t)}{\partial^2t^2}(1)  =\frac{-d(d-1)N+\sum_{i=1}^{N-1}(d+m_i)(d+m_i-1)}{2}
\nonumber\\
&=&\frac{-d(d-1)N+\sum_{i=1}^{N-1}(d(d-1)+(2d-1)m_i+m_i^2)}{2}\nonumber\\
&=&\frac{-d(d-1)N+d(d-1)(N-1)+(2d-1)\sum_{i=1}^{N-1}m_i+\sum_{i=1}^{N_1}m_i^2}{2}\nonumber\\
&=&\frac{-d(d-1)+(2d-1)d+\sum_{i=1}^{N-1}m_i^2}{2}\nonumber\\
&=& \frac{d^2+ \sum_{i=1}^{N-1}m_i^2}{2}.\nonumber
\end{eqnarray}
\qed

\begin{Corollary}\label{Bir_equi_linearsyz}
	Let $I\subset R=k[x,y,z]$ be a $d$-equigenerated homogeneous perfect ideal of height $2$ that is  generically complete intersection, but not a complete intersection. Then:
	\begin{enumerate}
		\item[{\rm (a)}] $\ell(I)=3.$ 
		\item[{\rm(b)}] The defining ideal of $\mathcal{F}(I)$ is generated in degrees $\geq 3$.
	\end{enumerate}
	Moreover, if the Rees algebra $\mathcal{R}(I)$ is Cohen-Macaulay then the following statements are equivalent:
	\begin{enumerate}
		\item[{\rm(i)}]  $I$ is linearly presented.
		\item[{\rm (ii)}] The rational map $\mathfrak{F}:\pp^2\to \pp^{N-1}$ defined by the linear  system $I_{d}$ is birational onto its image.
	\end{enumerate}
\end{Corollary}
\demo (a) This follows as usual from \cite{Cow-Nor}.

(b) This comes out of Corollary~\ref{Tchernev_Cor} since $I$ necessarily satisfies condition $G_3$. 

Now, we head out to the supplementary assertion.

(i)$\Rightarrow$(ii) By (a), this is a case of \cite[Theorem 3.2]{AHA}.

(ii) $\Rightarrow$ (i) Let $$0\to \bigoplus_{i=1}^{N-1} R(-d-m_i)\stackrel{\phi}\lar R(-d)^N\to I\to 0$$
be the minimal graded free resolution of $I$, where $m_i\geq 1$ for every $i$.  Since $I$ is generically a complete intersection, \cite[Theorem 6.6(b)]{Ram2} gives 
\begin{equation}\label{ram2_formula}
e(\mathcal{F(I)})\cdot \deg(\mathfrak{F}) = d^2-e(R/I).
\end{equation}
Applying the formula of Corollary~\ref{multiplicity_of_perfect_cod2} yields
\begin{eqnarray*}
	e(\mathcal{F(I)})\cdot \deg(\mathfrak{F}) &=& d^2-\frac{d^2+ \sum_{i=1}^{N-1}m_i^2}{2}\nonumber\\
	&=&\frac{d^2- \sum_{i=1}^{N-1}m_i^2}{2}\nonumber\\
	&=&\frac{\displaystyle\sum_{i=1}^{N-1}m_i^2+2\sum_{1\leq i< j\leq N-1}m_im_j- \sum_{i=1}^{N-1}m_i^2}{2}\nonumber\\
	&=&\sum_{1\leq i< j\leq N-1}m_im_j. 
\end{eqnarray*} 

By assumption, $\deg(\mathfrak{F})=1.$ Hence, 

$$e(\mathcal{F(I)})=\sum_{1\leq i< j\leq N-1}m_im_j$$
On the other hand, since $\mathcal{R}(I)$ is Cohen-Macaulay then, by items (a) and (b),  the hypotheses of Theorem~\ref{degreeRees_CM} are satisfied. Thus, $$e(\mathcal{F(I)})={N-1\choose 2}.$$ Hence,
$$\sum_{1\leq i< j\leq N-1}m_im_j={N-1\choose 2}.$$
This forces $ m_i=1$ for every $1\leq i\leq N-1$ as desired.
\qed

\begin{Remark}\rm
Among the hypotheses of Theorem~\ref{Bir_equi_linearsyz} used in the supplementary part, (a) can be dispensed with by instead adding to item (i)  the hypothesis that the entries of the syzygy matrix of $I$ generate the maximal ideal of $R$. This detour is afforded by \cite[Theorem 2.4 (a)]{lin_cod2}, of which the implication (ii) $\Rightarrow$ (i)  is a sort of weak converse. 
\end{Remark}

The following simple examples illustrate some aspects of the above results.

\begin{Example}\rm
	Let $I=(x^2,xy,xz,yz)\subset R=k[x,y,z]$.
	
(1) It is easy to see that $I$ is non-perfect since $x (x,y,z)\subset I$ and that $\mathcal{F}(I)$ is defined by an obvious degree $2$ binomial equation -- in particular, none of parts (b) and (c) of Theorem~\ref{degreeRees_CM} holds. 

Computation with \cite{M2} gives:


(2) $\mathcal{R}_R(I)$ is Cohen--Macaulay.

(3) $r(I)=1$, hence part (a) of Theorem~\ref{degreeRees_CM} does not hold either.
\end{Example}

\begin{Example}\rm
Let $I=(x^2-y^2,xy,xz,yz)\subset R=k[x,y,z]$.

(1) It is easy to see that $I$ is non-perfect (e.g., $x^2 (x,y,z)\subset I$) and a complete intersection at its unique minimal prime $(x,y)$. 

Computation with \cite{M2} gives:

(2) $I$ has maximal linear rank, but is not linearly presented; in particular, $I_2$ defines a birational map onto the image, thus showing that perfectness is essential for the implication (ii)$\Rightarrow$ (i) of Corollary~\ref{Bir_equi_linearsyz}. 

(3) $\mathcal{R}_R(I)$ is not Cohen--Macaulay, and yet $r(I)=2$.

(4)  $\mathcal{F}(I)$ is defined by a degree $3$ equation, hence both parts (b) and (c) of Theorem~\ref{degreeRees_CM} hold.
\end{Example}

\begin{Question}\rm
One wonders whether the following generalization of the supplementary part of Corollary~\ref{Bir_equi_linearsyz} holds: {\em under the same hypotheses, except that $I$ is not assumed to be generically a complete intersection, can one replace `linear presentation' by `sub-maximal linear rank'} ?
\end{Question}
Note that a perfect ideal of height $2$ has sub-maximal rank ($=N-2$) if and only if it is generated by linear syzygies plus one single syzygy of degree $\geq 2$.
The following simple example shows that the answer to the question is negative if the Rees algebra of $I$ is not assumed to be Cohen--Macaulay.

\begin{Example}\rm
Let $I$ be the ideal generated by the maximal minors of the following matrix$$\left[\begin{array}{cccc}
z^2&0&0&0\\
x^2&z^2&0&0\\
0&y^2&z&0\\
0&0&x&z\\
0&0&0&y
\end{array}\right]$$
\end{Example}
It is easy to see that $I$ is not generically a complete intersection (hence, does not satisfy $G_3$).
A calculation with \cite{M2} shows that $\mathcal{F}(I)$ is Cohen--Macaulay and generated in degree $3$.
This implies that the discussed rational map is birational onto the image, while the linear rank is only $2$. However, $\mathcal{R}_R(I)$ is not Cohen--Macaulay.



\noindent {\bf Addresses:}

\medskip

\noindent {\sc Zaqueu Ramos}\\
Departamento de Matem\'atica, CCET\\ 
Universidade Federal de Sergipe\\
49100-000 S\~ao Cristov\~ao, Sergipe, Brazil\\
{\em e-mail}: zaqueu@mat.ufs.br\\

\medskip

\noindent {\sc Aron Simis}\\
Departamento de Matem\'atica, CCEN\\ 
Universidade Federal de Pernambuco\\ 
50740-560 Recife, PE, Brazil\\
{\em e-mail}:  aron@dmat.ufpe.br

\end{document}